\documentclass{amsart}
\usepackage{math}
\usepackage{lineno,tikz,mathabx,mathrsfs,mathtools,hyperref}
\usetikzlibrary{cd,quotes,angles,decorations,arrows,calc,automata}
\usepackage[margin=42mm,a4paper]{geometry}

\newtheorem{mainthm}{Theorem}

\begin{document}

\title{Equations in wreath products}
\author{Leon Pernak}
\email{pernak@math.uni-sb.de}
\author{Ruiwen Dong}
\email{ruiwen.dong@uni-saarland.de}
\author{Jan Philipp W\"achter}
\email{j.ph.waechter@gmail.com}
\author{Laurent Bartholdi}
\email{laurent.bartholdi@gmail.com}
\address{Saarland University}
\date{October 1, 2024}
\thanks{The authors gratefully acknowledge support from the ERC AdG grant 101097307}
\begin{abstract}
  We survey solvability of equations in wreath products of groups, and prove that the quadratic diophantine problem is solvable in wreath products of Abelian groups. We consider the related question of determining commutator width, and prove that the quadratic diophantine problem is also solvable in Baumslag's finitely presented metabelian group. This text is a short version of an extensive article by the first-named authors.
\end{abstract}
\maketitle

\section{Introduction}
Let $G$ be a group. An \emph{equation} over $G$ is an expression $w(X_1,\dots,X_n)=1$ in which $w$ is a word in some variables $X_1,\dots,X_n$ and constants in $G$; for example, $X_1^2g^{-1}=1$ is an equation, which is solvable precisely when $g$ has a square root in $G$.

More formally, consider the free group $F_X$ over countably many generators $X_1,X_2,\dots$; then an equation over $G$ is given by an element $w\in F_X*G$, and this equation is solvable precisely when there exists a homomorphism $\sigma\colon F_X*G\to G$ that is the identity on $G$ and maps $w$ to $1$. This homomorphism $\sigma$ is called a \emph{solution} of the equation, identifying the solution with the values $(\sigma(X_1),\sigma(X_2),\dots)$.

\def\DP{{\mathcal{DP}}}
\def\QDP{{\mathcal{Q^+DP}}}
The fundamental question is then: \emph{For a group $G$, is there an algorithm that, given an equation in $G$, determines whether it admits a solution?} We call this problem $\DP_1$, the index indicating that a single equation is considered; more generally, the \emph{diophantine problem} $\DP$ in a group asks for an algorithm to determine solvability of a system of equations.  If $w$ is an equation involving no variable, i.e.\ a constant, then it is solvable if and only if $w$ represents the identity of $G$; thus a necessary condition for solvability of $\DP_1$ is that $G$ have solvable word problem. Some results:
\begin{itemize}
\item $\DP$ (and actually much more) is solvable in free groups, by the fundamental contributions of Makanin~\cite{makanin;equations-group} and Razborov~\cite{razborov;equations};
\item $\DP_1$ is solvable in the Heisenberg group of $3\times3$ upper-triangular unipotent integer matrices, by Duchin-Liang-Shapiro~\cite{duchin-liang-shapiro;equations};
\item $\DP_1$ is unsolvable in free metabelian groups, by Roman'kov~\cite{romankov;free-metabelian};
\item $\DP_1$ is unsolvable in free nilpotent groups of class $3$, by Truss~\cite{truss;equation};
\item $\DP$ is unsolvable in the wreath product $\Z\wr\Z$, by a recent result of Dong~\cite{dong;equations}.
\end{itemize}
(Wreath products will be recalled in the next section). Much current research attempts to delineate, among metabelian groups, those for which $\DP$ or $\DP_1$ is solvable.

An important subclass of equations has emerged, the \emph{quadratic equations}. They are those equations $w\in F_X*G$ such that every variable $X_i$ appears at most twice as a letter $X_i$ or $X_i^{-1}$ in $w$. Obviously if a variable appears only once then the equation is solvable, by solving for that variable, so no generality is lost in supposing that each of the variables $X_1,\dots,X_n$ appears exactly twice in $w$, and no other variable appears. If furthermore each of these $X_i$ appears once with exponent $+1$ and once with exponent $-1$, then $w$ is called an \emph{orientable quadratic equation}. The motivation comes from $2$-dimensional geometry: drawing $w$ on the boundary of a polygon and gluing each edge labeled $X_i$ to the one labeled $X_i^{-1}$ produces an orientable surface $\Sigma$ with elements of $G$ marked on boundary components. A solution is then nothing but a homomorphism $\pi_1(\Sigma)\to G$ which sends boundary components to specified values in $G$. The classification of surfaces implies that $w$ may be given the form
\begin{equation}\label{eq:quadratic}
  w = [X_1,X_2]\cdots[X_{2n-1},X_{2n}]c_1^{X_{2n+1}}\cdots c_p^{X_{2n+p}}
\end{equation}
where $\Sigma$ has genus $n$ and $p$ punctures, for some constants $c_1,\dots,c_p\in G$; here and below the commutator is $[g,h]=g^{-1}h^{-1}gh$. There is a parallel story for \emph{unorientable} quadratic equations, that we ignore for brevity. We denote by $\QDP$ the \emph{orientable quadratic diophantine problem}, that of determining solvability of an orientable quadratic equation (or equivalently an orientable quadratic system of equations, since the variables of different equations in the system are either disjoint or can be eliminated). Sample results include:
\begin{itemize}
\item $\QDP$ is solvable in the ``lamplighter group'' $C_2\wr\Z$ and in the Baumslag-Solitar groups $BS(1,k)=\langle a,t\mid a^t=a^k\rangle$, by Kharlampovich-L\'opez-Myasnikov~\cite{kharlampovich-lopez-myasnikov} (beware that the published version claims more than the revised, arXiv version);
\item $\QDP$ (and also the unorientable quadratic diophantine problem) is solvable in the Grigorchuk group, by Lysenok-Miasnikov-Ushakov~\cite{lysenok-miasnikov-ushakov;quadratic};
\item $\QDP$ is solvable in free metabelian groups, by Lysenok-Ushakov~\cite{lysenok-ushakov;quadratic-metabelian};
\item $\QDP$ is unsolvable in some (explicit, finitely generated) Abelian-by-cyclic group, again by Dong in the same article~\cite{dong;equations}.
\end{itemize}
(We shall not define the Grigorchuk group, but refer simply to that paper which contains a good introduction to it).

A closely related problem, in view of the commutator form of~\eqref{eq:quadratic}, is to determine the \emph{commutator width} of a group $G$, namely the minimal $n$ such that every element of $[G,G]$ is a product of at most $n$ commutators. Indeed if $n$ is bounded then $G$ admits ``commutator elimination'': the genus $g$ in~\eqref{eq:quadratic} may be assumed to be at most $n$; and conversely a clever argument in~\cite{lysenok-miasnikov-ushakov;quadratic} deduces for the Grigorchuk group that $n$ is finite, without expliciting its value (which turns out to be $3$, see~\cite{bartholdi-groth-lysenok;commutator}).

\subsection{Main results}
We recall that the wreath product $A\wr B$ of two groups is the semi-direct product $\prod'_B A\rtimes B$ of the restricted power (almost all coordinates are $1$) of $A$ by $B$, the action being by permutation of the coordinates. Elements of $A\wr B$ may be represented in the form $(f,b)$ with $f\colon B\to A$ and $b\in B$, and then the multiplication is given by
\[(f,b)\cdot(f',b')=(x\mapsto f(x)f'(x b),b b').\]
This basic construction in group theory is ubiquitous; for example, the group of symmetries of the cube is the wreath product $C_2\wr S_3$, the $2$-Sylow subgroup of the symmetric group $S_{2^n}$ is the iterated wreath product $C_2\wr\cdots\wr C_2$, etc.; one of the most-encountered examples is the ``lamplighter group'' $C_2\wr\Z$; and the Grigorchuk group may be defined via its imbedding $G\hookrightarrow G\wr C_2$.

\begin{mainthm}\label{thm:eqs}
  Let $A,B$ be two finitely generated Abelian groups. Then the quadratic diophantine problem $\QDP$ is solvable in $A\wr B$.
\end{mainthm}

\begin{mainthm}\label{thm:cw}
  Let $A,B$ be two finitely generated Abelian groups, with $A\neq1$. Then the commutator width of $A\wr B$ is $\lceil\rank(B)/2\rceil$.
\end{mainthm}

We concentrate on $A=\Z$ and $B=\Z^n$ in the sketches of the proofs; the difficulties to extend to the general case are more notational than technical, and definitely not conceptual. For the full treatment refer to the forthcoming article~\cite{dong-pernak;equations}.

We include, nevertheless, another example, the ``Baumslag group'', to show the extent of our method, the difficulties that may arise to treat metabelian groups beyond wreath products, and hopefully to convince the reader that ``smaller'' groups may have harder diophantine problems than ``larger'' or ``freeer'' ones. The Baumslag group admits as presentation
\begin{equation}\label{eq:baumslag}
  \Gamma\coloneq\langle a,t,u\mid [t,u], [a,a^t], u^u = a a^t\rangle
\end{equation}
and is remarkable as an example of a finitely presented metabelian group containing the wreath product $\Z\wr\Z$ (as the subgroup $\langle a,t\rangle$). We show:
\begin{mainthm}\label{thm:baumslag}
  $\QDP$ is solvable in the Baumslag group $\Gamma$, and $\Gamma$ has commutator width $1$.
\end{mainthm}

\section{Proofs}
For an Abelian group $A$ that has the structure of a commutative ring, the wreath product $A\wr B$ has a natural expression in terms of the group ring: it is the semidirect product $A[B]\rtimes B$ of $B$'s group ring with $B$, the action of $B$ in its group ring being given by right multiplication. To simplify the notation, we restrict ourselves to $A=\Z$.

For every quotient map $\pi\colon B\twoheadrightarrow Q$, there is a natural ideal of the group ring
\[\varpi_Q\coloneq\ker(\Z[B]\twoheadrightarrow\Z[Q]).\]
If furthermore $B$ is free Abelian, say $B=\Z^d$, then its group ring $\Z[B]$ is isomorphic to the ring of Laurent polynomials $\Z[Z_1^{\pm1},\dots,Z_d^{\pm1}]$. The group $B$ is identified in $\Z[B]$ with Laurent monomials, and we write $Z^g$ for the monomial corresponding to $g\in B$; so $\varpi_{B/N}$ is generated as an ideal by $\{1-Z^h:h\in N\}$.

\subsection{Spherical equations}
These are the quadratic equations with genus $0$, in other words of the form $w=\prod_{i=1}^p g_i^{X_i}$ for some elements $g_i\in\Z\wr B$. Writing $g_i=(f_i,t_i)$ with $f_i\in\Z[B],t_i\in B$, we observe that, if $t_i=1$, then the corresponding variable $X_i$ may be assumed to belong to $B$ since the component of $X_i$ in $\Z[B]$ commutes with $f_i$.

Evidently the constants in $w$ must satisfy $\prod t_i=1$ if there is any chance of a solution. For any $j\in\{1,\dots,p\}$, write $\overline B=B/\langle t_j\rangle$; we claim that $w$ is solvable if and only if its projection to $\Z\wr\overline B=\Z[B]/\varpi_{\overline B}\rtimes\overline B$ is solvable. In the non-trivial direction, let $\overline\sigma$ be a solution in $\Z\wr\overline B$; by the remark above $\overline\sigma(Z_j)$ may be assumed to belong to $\overline B$. Lift $\overline\sigma$ to $\sigma\colon F_X\to\Z\wr B$, write $\sigma(Z_j)\eqcolon u\in B$, and write $(f',u'),(f'',u'')$ respectively for the subproducts in $w$ on the left and right of $(f_j,t_j)^u$ in $\sigma(w)$; we obtain
\[f\coloneq f'+f_j\cdot u'+f''\cdot(u' u)\in\varpi_{\overline B}=\langle1-Z^{t_j}\rangle\]
and we adjust $\sigma(Z_j)\coloneq (-f/(1-Z^{t_j}),u)$ to obtain a solution. Repeating, we reduce to equations for which the constants are all of the form $g_i=f_i\in\Z[B]$.

We are now reduced to the following problem in group rings: \emph{given $f_1,\dots,f_p\in\Z[B]$, determine whether there exist translates $f_1\cdot u_1,\dots,f_p\cdot u_p$ which sum to $0$}. This can only happen if the supports of the $f_i\cdot u_i$ overlap in independent clusters summing to $0$, so we may bound the lengths of the $u_i$ in terms of the $f_i$'s supports. Spherical equations are thus decidable.

\subsection{Commutators in wreath products}
An easy calculation gives, for $(f,t),(f',t')\in\Z\wr B$, the commutator
$[(f,t),(f',t')] = ((1-Z^{t'})f-(1-Z^t)f',1)$,
so the set of commutators in $\Z\wr B$ is
\[\bigcup_{t,t'\in B}\langle1-Z^t,1-Z^{t'}\rangle\times\{1\},\]
and the set of $n$-fold products of commutators is
\[\bigcup_{t_1,\dots,t_{2n}\in B}\langle1-Z^{t_1},\dots,1-Z^{t_{2n}}\rangle\times\{1\}.\]

We deduce that the commutator width of $\Z\wr B$ is $\lceil\rank(B)/2\rceil$, proving Theorem~\ref{thm:cw}: in one direction, products of $\lceil\rank(B)/2\rceil$ commutators contain $\varpi_1\times\{1\}=\gamma_2(\Z\wr B)$, and in the other direction the element $d-Z_1-\cdots-Z_d$ does not belong to any ideal generated by less than $d$ terms of the form $1-Z^t$.

We also deduce that an equation of the form $w=[X_1,Y_1]\cdots[X_n,Y_n]w'$ in $\Z\wr B$ has a solution if and only if there exists a lattice $L\le B$ of rank at most $2n$ such that $w'$ has a solution in $\Z\wr(B/L)$. Using the same reductions as above, we are again led to a group ring problem: \emph{given $f_1,\dots,f_p\in\Z[B]$, determine whether there exist translates $f_1\cdot u_1,\dots,f_p\cdot u_p$ which sum to an element of $\varpi_{B/L}$}. Again this can only happen if the supports of the $f_i\cdot u_i$ overlap in clusters, leading to bounds on the lengths of the $u_i$ and only finitely many cases need be considered; and the rank of the space spanned by the corresponding shifted supports is computable, for example by evaluating appropriate determinants (the Pl\"ucker relations). We have completed the proof of Theorem~\ref{thm:eqs}.

\section{The Baumslag group}
In~\cite{baumslag;finitely-presented-metabelian}, Baumslag gives the group $\Gamma$ in~\eqref{eq:baumslag} as an example of a finitely presented metabelian group containing an infinitely presented subgroup. He later shows that every finitely generated metabelian group can be imbedded in a finitely presented metabelian group~\cite{baumslag;metabelian-embedding}.

Here $\Gamma$ is the quotient of $\Z\wr\Z^2$ by the relations $f+f\cdot t+f\cdot u=0$ for every $f\in\Z[\Z^2]$, namely it may be written as
\[\Gamma=\Z[Z_1^{\pm1},Z_2^{\pm1}]/(1+Z_1-Z_2)\rtimes\Z^2.\]
On the one hand, since $\Gamma$ is a quotient of $\Z\wr\Z^2$, its commutator width is at most $1$, so is precisely $1$ because it is not Abelian. On the other hand, to solve equations we follow the arguments above in the ring $\Z[Y^{\pm1}]$ acted upon by $Z_1=Y$ and $Z_2=Y+1$. In the case of genus $>0$, the whole derived subgroup $\gamma_2(\Gamma)$ is expressible with a single commutator, so the equation to be solved actually takes place in $\Gamma/\gamma_2(\Gamma)\cong\Z^3$; therefore the interesting case is that of a spherical equation.

By the reductions from the previous section, the spherical equation may be assumed to have constants in $K\coloneq\langle a\rangle^\Gamma=\ker(\Gamma\twoheadrightarrow\Z^2)$, and the variables may be assumed to take values in $\langle t,u\rangle$, so we arrive at the following problem: \emph{given $f_1,\dots,f_p\in K$, determine whether there exist translates $f_1\cdot u_1,\dots,f_p\cdot u_p$ with $u_i\in\langle t,u\rangle$ which sum to $0$}. Since $K=\Z[Z_1^{\pm1},Z_2^{\pm1}]/(1+Z_1-Z_2)$ is isomorphic to the ring of Laurent polynomials $\Z[Y^{\pm1}]$, the problem becomes: \emph{Given $f_1,\dots,f_p\in\Z[Y^{\pm1}]$, determine whether there exist $a_1,b_1,\dots,a_p,b_p\in\Z$ with $Y^{a_1} (Y+1)^{b_1} f_1 + \dots + Y^{a_p} (Y+1)^{b_p} f_p = 0$.}

Without loss of generality, by considering all possible orderings of the $f_i$, we may suppose $b_1 \ge b_2 \ge \dots \ge b_p = 0$. Suppose by induction on $k$, starting at $k=p$ downwards, that we have an upper bound on $b_{k+1}, b_{k+2},\dots, b_p$. Then $b_k$ can be bounded in the following way. Since $b_{k+1}, b_{k+2},\dots,b_p$ are bounded from above, the number of monomials appearing in the polynomial $\Phi\coloneq Y^{a_{k+1}} (Y+1)^{b_{k+1}} f_{k+1} + \dots + Y^{a_p} (Y+1)^{b_p} f_p$ is also bounded, say by $B$. Then since $Y^{a_1} (Y+1)^{b_1} f_1 + \dots + Y^{a_k} (Y+1)^{b_k} f_k = -\Phi$, we see that $(Y+1)^{b_k}$ divides $\Phi$. If $\Phi=0$ then we have two shorter equations to solve, and recurse; otherwise, $b_k \le B-1$ by~\cite[Lemma~4.1]{giesbrecht-roche-tilak;sparse}, and the induction on $k$ proceeds. It follows that the $b_i$'s are all bounded, and trying them all yields equations of the form $Y^{a_1} \widetilde f_1+\dots+Y^{a_p}\widetilde f_p=0$, that we have already considered and know how to solve. Theorem~\ref{thm:baumslag} is proven.

It is the nature of the quotient $\Z[Z_1^{\pm1},Z_2^{\pm1}]/(1+Z_1-Z_2)\cong\Z[Y^{\pm1}]$ that makes the problem solvable; a precise formulation of an algebraico-geometric condition guaranteeing this is still wanting.

\bibliographystyle{amsalpha}
\bibliography{equations}

\end{document}